\documentclass{birkjour}
\usepackage[T1]{fontenc} 
\usepackage{latexsym} 
  \usepackage[all]{xy}
  \usepackage{amsfonts} 
  \usepackage{amsthm} 
  \usepackage{amsmath} 
  \usepackage{amssymb}
  \usepackage{pifont}  
  \usepackage{enumerate}
  
 \newtheorem{thm}{Theorem}[section]

 \theoremstyle{definition}
 
 \theoremstyle{remark}

 \numberwithin{equation}{section}
 
 \def\bs{\begin{statement}}
\def\es{\end{statement}}
  \swapnumbers
  \newtheorem{statement}[thm]{}

  \newcounter{zlist}

  \newcounter{blist}

  \newcounter{rlist}

\def\AA{{\mathbb A}}

\def\CC{{\mathbb C}}

\def\NN{{\mathbb N}}

\def\ZZ{{\mathbb Z}}

\def\tt{{\mathfrak T}}

\newcommand{\Cc}{\mathcal{C}}

\def\O{{\bf O}}

\def\*C{{}^*\hspace*{-1pt}{\Cc}}

\def\text#1{{\rm {\rm #1}}}

 \def\1{\mathbf{1}}

 \def\supp{\mathrm{supp}}

\begin{document}

%
%
%
%
%
%
%
%
%

\title[Skew derivations on down-up algebras]
 {Skew derivations on 
 down-up algebras}

\author{Munerah Almulhem}

\address{%
 Department of Mathematics, Swansea University, 
  Swansea SA2 8PP, U.K. \ \newline
 Department of Mathematics, Imam Abdulrahman bin Faisal University, 
  Dammam 34212, K.S.A.} 
  \email{844404@swansea.ac.uk} 
  
  \author{Tomasz Brzezi\'nski}

\address{
Department of Mathematics, Swansea University, 
  Swansea SA2 8PP, U.K.\ \newline \noindent
Department of Mathematics, University of Bia{\l}ystok, K.\ Cio{\l}kowskiego  1M,
15-245 Bia\-{\l}ys\-tok, Poland}

\email{T.Brzezinski@swansea.ac.uk}

\subjclass{Primary 16S38; Secoundary 16W25; 58B32}

\keywords{generalized down-up algebra; generalized Weyl algebra; skew derivation}

\date{October 2016}

\begin{abstract}
A class of skew derivations on complex Noetherian generalized down-up algebras $L=L(f,r,s,\gamma)$ is constructed.
\end{abstract}

\maketitle

\section{Motivation and introduction}
Determining which classes of non-commutative algebras correspond to smooth non-commutative varieties or manifolds is one of the outstanding problems of non-commutative geometry or NCG. In a recent proposal \cite{BrzSit:smo} it is argued that a possibility of constructing a suitable graded differential algebra or a differential structure should determine smoothness of a non-commutative variety. This kind of smoothness is referred to as {\em differential}. The proposal involves a strict Poincar\'e type duality between differential and integral forms (the reader might like to consult \cite{Brz:lec} for a concise explanation of terms involved) and is {\em constructive} in nature. Despite some recent progress in  uncovering functorial ways of checking differential smoothness \cite{BrzLom:dif} one needs to study algebras on case by case basis and construct suitable differential and integral complexes. In the background for such complexes one often finds {\em skew derivations}, which can be understood as vector fields in NCG. 

Motivated by the role skew derivations play in NCG, in particular in constructing smooth structures, in \cite{AlmBrz:ske} we undertook a detailed study of skew derivations on generalized Weyl algebras \cite{Bav:gen}. Despite their almost na\"ive simplicity, generalized Weyl algebras include an astounding number of examples that appeared and still appear in NCG, and in other algebraic contexts. In this note we concentrate on one such example in the latter class.

Generalized down-up algebras were introduced by Cassidy and Shelton in \cite{CasShe:bas} as a generalization of the down-up algebras of Benkart and Roby  \cite{BenRob:dow}. Significant, classical examples of generalized down-up algebras include the enveloping algebra of $sl_2$, and the enveloping algebra of the $3$-dimensional Heisenberg Lie algebra. More recent examples are various algebras similar to the enveloping algebra of $sl_{2}$ such as those introduces by Smith \cite{Smi:clas}, Witten \cite{Wit:gau},  Le Bruyn \cite{LeB:con} and Rueda \cite{Rue:som}. Generalized down-algebras are affine algebras of Gelfand-Kirillov dimension three, i.e.\ they can be interpreted as coordinate algebras of three-dimensional non-commutative varieties, and in the case in which they are Noetherian domains (on which we focus in this note), it is natural to ask whether they are differentially smooth. In this note we do not attempt to answer this question, but rather provide the first paving stones for a path that might lead to an answer by constructing a class of skew derivations.

\section{Generalized down-up and Weyl algebras}\label{sec.Noeth}
Let $\Bbbk$ be an algebraically closed field of characteristic zero. 
Fix scalars $r,s,\gamma \in\Bbbk$ and a polynomial $f \in \Bbbk[X]$. A {\em generalized down-up algebra} $L:=L(f,r,s,\gamma)$  \cite{CasShe:bas} is a unital associative $\Bbbk$-algebra generated by $d$, $u$ and $h$, subject to relations:
\begin{equation*}
dh-rhd+\gamma d=0, \quad
hu-ruh+\gamma u=0, \quad
du-sud+f(h)=0.
\end{equation*}
When $f$ has degree one, we retrieve all down-up algebras of \cite{BenRob:dow}. 
It was explained in \cite{CasShe:bas} that $L$ has Gelfand-Kirillov dimension three and it is Noetherian if and only if it is a domain, if and only if $ rs\neq 0$. From now on, we always assume $rs \neq 0$, and we also set $\Bbbk = \CC$. 

$L$ is said to be {\em conformal} if there exists a polynomial $g\in\CC[X]$ such that $f(X)=sg(X)-g(rX-\gamma)$, e.g.\  $L(0,r,s,\gamma)$ is conformal. By \cite[lemma~ 2.8]{CasShe:bas}, $L$ is conformal if $s \neq r^{i}$ for all $0\leq i \leq \deg(f)$ (this is a sufficient, but not necessary condition for conformality). If $r\neq 1$, then, by \cite[Propostion~1.7]{CarLop:aut},  any generalized down-up algebra is isomorphic to one with $\gamma =0$. A generalized down-up algebra $L(f,r,s,0)$ is conformal if and only if $r\neq s^i$ for all $i$ in the support of  $f$ (i.e.\ for those $i\in \NN$, for which $X^i$ has a non-zero coefficient in the expansion of $f(X)$). In this case the support of $g$ can be assumed to be equal to that of $f$, and then $g$ is uniquely determined by $f$; in particular $\deg(f) = \deg(g)$.

Given an algebra  $R$, an automorphism $\varphi$ of $R$ and a central element $a\in R$, the {\em generalized Weyl algebra} $R(a,\varphi)$ is a ring extension of $R$ generated by $x$ and $y$, subject to the relstions:
 \begin{equation}\label{gW}
 xy = \varphi(a), \quad yx = a, \quad xr = \varphi(r)x, \quad yr = \varphi^{-1}(r)y, 
 \end{equation}
Generalized Weyl algebras were introduced and studied by Bavula in \cite{Bav:gen}. Any $R(a,\varphi)$ is a $\ZZ$-graded algebra with $R$ contained in the degree-zero part and $\deg(x) = 1$, $\deg(y) =-1$. If $R$ is a Noetherian algebra which is a domain  and $a\neq 0$ then $R(a,\varphi)$ is a Noetherian domain. Noetherian generalized down-up algebras can be presented as generalized Weyl algebras as follows. Set $a=ud$, let $R$ be the commutative polynomial algebra $\CC[h,a]$ and define the automorphism $\varphi$ by the rules $\varphi(h)=rh-\gamma$ and $\varphi(a)=sa-f(h)$. Then 
$$
\CC[h,a]\left(a,\varphi\right)\cong L(f,r,s,\gamma).
$$

Henceforth we assume that $r$ is not a root of unity, $rs\neq 0$ and $s\neq r^i$, for all $i\in \NN$. Thus $L$ is Noetherian and conformal, and with no loss of generality (up to isomorphism) we can also assume that $\gamma=0$. In this case there is a more convenient presentation of $L$ as a generalized Weyl algebra. Let $g$ be the unique polynomial (of the same degree and with the same support as $f$) such that $f(X)=sg(X)-g(rX)$.
Let $a=ud$ and $k=a-g(h)$. Then $\CC[a,k] = \CC[h,k]$ and the automorphism $\varphi$ acts on $k$ by 
\begin{eqnarray*}
\varphi(k)=\varphi(a-g(h))=sa-f(h)-g(rh)=sa-sg(h)=sk.
\end{eqnarray*}
Therefore, $L$ is presented as the generalized Weyl algebra $R(\varphi,k+g(h))$, where $R=\CC[h,k]$ and $\varphi$ is the automorphism of $R$ defined by 
\begin{equation}\label{phi}
\varphi(h)=rh, \qquad \varphi(k)=sk.
\end{equation}
 The relations are thus:
 \begin{subequations}\label{rel.skew}
\begin{equation}
 xy=sk+g(rh),  \qquad yx=k+g(h),
 \end{equation}
 \begin{equation}  
 xp(h,k) = p(rh,sk)x, \qquad yp(h,k) = p(r^{-1}h,s^{-1}k)y, 
\end{equation}
\end{subequations}
for all $p(h,k)\in \CC[h,k]$. 
$L(f,r,s,0)$ is recovered from $R(k+g(h),\varphi)$  by the isomorphism $h\mapsto h$, $k\mapsto ud-g(h)$, $x\mapsto d$, $y\mapsto u$.
 
\section{Skew derivations on generalized down-up algebras}\label{sec.skew}
For any algebra  $A$, a {\em (right) skew derivation} or a {\em $\sigma$-derivation}  is a pair $(\partial, \sigma)$ consisting of an algebra endomorphism $\sigma: A\to A$ and a linear map $\partial: A\to A$ that satisfies the $\sigma$-twisted Leibniz rule, for all $a,b\in A$,
\begin{equation}
\partial(ab) = \partial(a)\sigma(b) + a\partial(b).
\end{equation}
A skew-derivation is said to be {\em inner} if it is given by a twisted commutator with an elemnent of $A$, i.e.\ for all $a\in A$, $\partial(a) = b\sigma(a)- ab$.

In \cite[Theorem~3.1]{AlmBrz:ske} we constructed a large class of skew derivations on arbitrary generalized Weyl algebras. The aim of this section is apply this construction to Notherian conformal generalized down-up algebras.

Let $A = R(a,\varphi)$ be a generalized Weyl algebra. Any algebra authomorphism $\sigma$ of $R$ such that 
\begin{equation}\label{sig.phi}
\sigma\circ\varphi = \varphi\circ\sigma \quad \mbox{and} \quad  \sigma(a)=a,
\end{equation}
 can be extended to an automorphism $\sigma_\mu$ of $A$ by setting,
\begin{equation}\label{deg.ext}
\sigma_\mu\mid_R = \sigma, \qquad \sigma_\mu(x) = \mu^{-1} x, \qquad \sigma_\mu(y) = \mu y,
\end{equation}
where $\mu\in \CC^\times$. The automorphism $\sigma_\mu$ is called a {\em degree-counting extension} of $\sigma$ of {\em coarseness} $\mu$.

In  \cite{AlmBrz:ske}, in addition to general inner skew derivations, two classes of {\em elementary} $\sigma_\mu$-skew derivation have been identified: 
\begin{description}
\item[(1) $c$-type derivations] To any $w\in \ZZ$ and $c_w\in R$, such that
\begin{equation}\label{c.phi}
bc_w = c_w\varphi^w(\sigma(b)), \qquad \mbox{for all $b\in R$},
\end{equation}
one can associate a (unique if $R$ is a domain) derivation $\partial_n^c$ on $A$ of $\ZZ$-degree or {\em weight} $w$ and such that $\partial_w^c\left(R\right) = 0$.
\item[(2) $\alpha$-type derivations] For any $w\in \ZZ\setminus\{0\}$, any $\varphi^w\circ\sigma$-skew derivation $\alpha_w$ of $R$ such that
\begin{equation}\label{alpha.phi}
\alpha_w\circ \varphi = \mu\, \varphi\circ \alpha_w,
\end{equation}
can be (uniquely) extended to a $\sigma_\mu$-skew derivation $\partial_w^\alpha$ of $A$ of $\ZZ$-degree $w$ such that $\partial_w^\alpha(x)=0$ if $w>0$ and $\partial_w^\alpha(y)=0$ if $w<0$. Also a $\sigma$-skew derivation $\alpha_0$ satisfying \eqref{alpha.phi} and such that $\alpha_0(a) = ca$, for some $c\in R$ such that $bc = c\sigma(b)$, for all  $b\in R$, can be extended to a $\sigma_\mu$-skew derivation of $A$ of the $\ZZ$-degree or weight zero. 
\end{description}

We now specify to conformal Noetherian generalized down-up algebras $L=L(f,r,s,0)$ viewed as generalized Weyl algebras $\CC[h,k](\varphi,k+g(h))$. Our standing assumptions are that the parameters $r,s \in \CC$ are neither zero nor roots of unity and  $s \neq r^{i}$, for all $ i \in \NN$. 

The identity map is the only automorphism of  $\CC[h,k]$ that satisfies \eqref{sig.phi}. Its extension to the whole of $L$ is thus given by:
\begin{eqnarray}\label{sig.mu}
\sigma_\mu(h)=h, \qquad \sigma_\mu(k)=k, \qquad \sigma_{\mu}(x)=\mu^{-1}x, \qquad \sigma_{\mu}(y)=\mu y.
\end{eqnarray}  
This automorphism is non-trivial if and only if $\mu \neq1$, which we assume from now on. 

Since $\sigma$ is the identity map,  the $w$-fold self-composition of $\varphi$ \eqref{phi}, $\varphi^w$, sends  a polynomial $p(h,k)$ to $p(r^wh,s^wk)$, and $\CC[h,k]$ is commutative, the equations \eqref{c.phi} can only be satisfied non-trivially if $w=0$. Consequently, there are only weight zero $\sigma_\mu$-skew $c$-type derivations of $L$: 
\begin{equation}\label{c0}
\partial^{c}_0(x) =c_{0}(h,k)x, \qquad  \partial^{c}_0(y) = -\mu c_{0}(r^{-1}h,s^{-1}k)   y ,
 \end{equation}
 for any $c_{0}(h,k)\in \CC[h,k]$; see \cite[Lemma~3.3]{AlmBrz:ske}. This derivation is inner, provided there exists $p(h,k) \in \CC[h,k]$ such that 
\begin{equation}\label{c.inner} 
 c_{0}(h,k)=\mu^{-1} p(h,k)-p(rh,sk).
 \end{equation}
By comparing the coefficients at powers of $h$ and $k$, one easily finds that \eqref{c.inner} can be solved if, and only if, 
$$
\mu^{-1}\neq r^\beta s^\gamma, \qquad \mbox{for all $(\beta,\gamma) \in \supp(c_0(h,k))$}.
$$
Therefore, if $\mu = r^{-\beta}s^{-\gamma}$ for some $\beta,\gamma \in \NN$, then $c_0(h,k) \sim h^\beta k^\gamma$  gives rise to a non-inner skew derivation \eqref{c0} on $L$ that vanishes on $\CC[h,k]$.

To construct $\alpha$-type derivations we need to consider skew derivations $\alpha_w$ of $\CC[h,k]$   twisted by $\varphi^w$. Any such skew derivation is fully determined by its values on generators of $\CC[h,k]$, i.e.\ there exists a unique $\varphi^w$-skew derivation $\alpha_w$ which on $h$ and $k$ is equal to any chosen pair of elements of $\CC[h,k]$, say
\begin{equation}\label{alpha.h.k}
\alpha_w(h)= \sum_{i,j}
\alpha^w_{h,i\,j} \,\,h^{i}k^{j},  \qquad 
\alpha_w(k)=\sum_{m,n}
\alpha^w_{k,m\,n}\,\, h^{m}k^{n}.
\end{equation} 
The value of $\alpha_w$ on any polynomial is then obtained by applying the $\varphi^w$-twisted Leibniz rule. The resulting skew derivation $\alpha_w$ satisfies condition \eqref{alpha.phi} if and only if this condition is satisfied for $\alpha_w$ evaluated on generators of $\CC[h,k]$. This in turn is equivalent to equations,
\begin{subequations}\label{alpha.phi.d-u}
\begin{equation}
r \sum_{i,j} \alpha^w_{h,i\,j} \,\, h^{i}k^{j} = \mu \sum_{i,j} r^{i}s^{j}\alpha^w_{h,i\,j} \,\, h^{i}k^{j},
\end{equation}
\begin{equation} 
 s\sum_{m,n} \alpha^w_{k,m\,n} \,\, h^{m}k^{n} = \mu \sum_{m,n} r^{m}s^{n}\alpha^w_{k,m\,n} \,\, h^{m}k^{n}. 
\end{equation}
\end{subequations}
Equations \eqref{alpha.phi.d-u} yield the following constraints:
\begin{equation}\label{con.r.s}
r^{i}s^{j}=\dfrac{r}{\mu},\; (i,j) \in \mathrm{supp}(\alpha_n(h)); \;\;\;
r^{m} s^{n}=\dfrac{s}{\mu}, \; (m,n) \in  \mathrm{supp}(\alpha_n(k)),
\end{equation}
where, for any polynomial $\pi\in \CC[h,k]$, $\supp(\pi)$ denotes its {\em support}, defined as the set of all pairs $(m,n)\in \NN$, for which $\pi$ has a non-trivial $h^mk^n$-term.
Constraints \eqref{con.r.s} restrict the supports of $\alpha_w(h)$ and $\alpha_w(k)$ as well as possible values of $\mu$. Set
\begin{equation}\label{r.s}
r = s^{b_1}, \qquad \mu^{-1}=s^{b_{2}}, \qquad \mbox{where $b_{1} \in \CC^\times \setminus \{\frac{1}{q}\; |\; q\in \NN^\times\}$, $b_{2} \in \CC^\times$}.
\end{equation}
The removal of zeros and the fractions $1/q$ ensures that $s$ is not a natural power of $r$, which is one of our standing assumptions (needed for $L$ to be both conformal and Noetherian), and that $\mu\neq 1$.  Note that the vector $\mathbf{b}:=(b_1,b_2)$ is not determined uniquely by $r,s$ and $\mu$. Solving constraints \eqref{con.r.s} we obtain the  following forms of $\alpha_w(h)$ and $\alpha_w(k)$:
\begin{equation}\label{alpha.all}
\alpha_w(h)=\sum_{i\in I} \alpha^w_{h,i} h^{i}k^{b_{2}+(1-i)b_{1}}, \quad 
\alpha_w(k)=\sum_{m\in J} \alpha^w_{k,m} h^{m} k^{b_{2}-mb_{1}+1},
\end{equation} 
where $I$ and $J$ are any finite subsets of 
\begin{equation}\label{indices}
I_{\mathbf{b}} := \{t\in \NN\; |\; b_{2}+(1-t)b_{1}\in \NN\}, \quad J_{\mathbf{b}} := \{t\in \NN\; |\; b_{2}-tb_{1}+1\in \NN\},
\end{equation}
respectively. What $I_{\mathbf{b}}$ and $J_{\mathbf{b}}$ are depends on $\mathbf{b}= (b_1,b_2)$. For example if $b_1$ is a negative integer and $b_2$ is a positive integer, then  $J_{\mathbf{b}}=\NN$, while $I_{\mathbf{b}}=\NN$ if $b_2 \geq -b_1$ and $I_{\mathbf{b}}=\NN\setminus \{0\}$ otherwise.  If $b_1$ is a positive integer or if $b_2$ is negative, then the indexing sets $I_{\mathbf{b}}$ or $J_{\mathbf{b}}$ could become finite or even empty. For example, if 
both $b_1$ and $b_2$ are positive, then
$$
I_{\mathbf{b}} \!= \!
\begin{cases} 
\{0,1\} & b_1>b_2,\cr
\{0, 1,2\} & b_1=b_2,\cr
\{0,\dots , q\!+\!1\} & b_1<b_2, 
\end{cases}
\; \;\;
J_{\mathbf{b}} \!= \! 
\begin{cases} 
\{0\} & b_1>b_2+1,\cr
\{0, 1\} & b_1=b_2, b_2+1,\cr
\{0,\dots , q\!+\!\delta_{\rho, b_1-1}\} & b_1<b_2,
\end{cases}
$$
where $q$ is the quotient and $\rho$ the remainder of the division of $b_2$ by $b_1$. 

If $b_1$ is positive and $b_2$ is negative such that $b_1<|b_2|$, then $I_{\mathbf{b}} = \emptyset$, hence there are no solutions to \eqref{con.r.s}.

The sets $I_{\mathbf{b}}$ and $J_{\mathbf{b}}$ may be non-empty even if the components of $\mathbf{b}$ are not integral. For example, if $b_1=b_2$, then $2\in I_{\mathbf{b}}$ and $1\in J_{\mathbf{b}}$.

The formulae \eqref{alpha.all} and \eqref{indices} lead to all possible elementary $\alpha$-type $\sigma_\mu$-skew derivations of non-zero weight on $L$, which in view of  \cite[Lemma~3.2]{AlmBrz:ske} take the following form:
\begin{subequations}\label{partial.a}
\begin{equation}
\partial^\alpha_w (p(h,k)) = 
\begin{cases}
\alpha_w(p(h,k)) x^w & w>0,\\
\alpha_w(p(h,k)) y^{-w} & w<0,
\end{cases}
\end{equation}
\begin{equation}
\partial^\alpha_w (x) = 
\begin{cases}
0  & w>0,\\
\mu^{-1}\alpha_w(sk+g(rh)) y^{-w-1} & w<0,
\end{cases}
\end{equation}
\begin{equation}
\partial^\alpha_w (y) = 
\begin{cases}
\mu\alpha_w(k+g(h)) x^{w-1} & w>0,\\
0 & w<0,
\end{cases}
\end{equation}
\end{subequations}
for all $p(h,k)\in \CC[h,k]$. Linear combinations of skew derivations, \label{c0} and \eqref{partial.a} contain all non-inner derivations of a Noetherian conformal down-up algebra $L(f,r,s,0)$ twisted by the automorphism \eqref{sig.mu} and vanishing on $x$ or $y$ or both $h$ and $k$.

In the weight-zero case in addition to \eqref{con.r.s}, one also needs to require that $\alpha_0(k + g(h))$ contains factor $k + g(h)$. In contrast to the non-zero weight case the values of $\alpha_0$ on $h$ and $k$ depend heavily on the choice of $g(h)$, and we take this as an excuse for not including them here.

\subsection*{Acknowledgments}
The research of the second author is partially supported by the Polish National Science Centre grant 2016/21/B/ST1/02438.

\end{document}